\documentclass[12pt]{article}
\usepackage{amssymb,fullpage,pstricks,graphicx,times}
\newtheorem{theorem}{Theorem}[section]

\newtheorem{lemma}[theorem]{Lemma}

\renewcommand{\P}{\mbox{$\mathbb{P}$}}
\newcommand{\E}{\mbox{$\mathbb{E}$}}

\def\phib{{\overline{\Phi}}}
\def\phibi{{\overline{\Phi}^{-1}}}
\def\ta{{T^{(1)}_j }}
\def\tb{{T^{(2)}_j }}

\let\hat\widehat
\let\tilde\widetilde

\begin{document}
\baselineskip 12pt

\begin{center}
{\Large{\bf Weighted Hypothesis Testing}}
\end{center}
\begin{center}
  {\sf Larry Wasserman and Kathryn Roeder}\footnote{Research supported
    by National Institute of Mental Health grants MH057881,MH066278,
    MH06329 and NSF Grant AST 0434343. The authors thank Jamie Robins
    for
    helping us to clarify several issues.}\\
  {\sf Carnegie Mellon University}\\
\end{center}
\begin{center}
{\sf{\Large{\gray April 7, 2006}}}
\end{center}

\begin{quote}
The power of multiple testing procedures can be increased by using
weighted p-values (Genovese, Roeder and Wasserman 2005).
We derive the optimal weights and
we show that the power is remarkably robust to misspecification
of these weights.
We consider two methods for choosing weights in practice.
The first, external weighting, is based on prior information.
The second, estimated weighting, uses the data to choose weights.
\end{quote}

\baselineskip 18pt

\section{Introduction}

The power of multiple testing procedures can be increased by using
weighted p-values (Genovese, Roeder and Wasserman 2005).
Dividing each p-value $P$ by a weight $w$ 
increases the probability of rejecting some hypotheses.  
Provided the weights have mean one, familywise error control
methods and false discovery control methods maintain their frequentist
error control guarantees.

The first such weighting scheme appears to be Holm (1979).  Related
ideas are in Benjamini and Hochberg (1997) and Chen et al (2000).
There are, of course, other ways to improve power
aside from weighting.  Some notable recent approaches include
Rubin, van der Laan and Dudoit (2005),
Storey (2005), Donoho and Jin (2004) and Signoravitch (2006).
Of these, our approach is closest to
Rubin, van der Laan and Dudoit (2005), hereafter, RVD.
In fact, the optimal weights derived here,
if re-expressed as cutoffs for test statistics,
turn out to be identical to the cutoffs
derived in RVD.
Our main contributions beyond RVD
are (i) a careful study of potential power losses
due to departures from the optimal weights, 
(ii) robustness properties of weighted methods, and
(iii) recovering power after using data splitting to estimate the weights.
An important distinction between
this paper and RVD
versus Storey (2005) is that Storey uses a slightly different
loss function and he requires a common cutoff for all test statistics.
This allows him to make an elegant connection with the Neyman-Pearson lemma.
In particular, his method automatically adapts from
one-sided testing to two-sided testing depending on the
configuration of means.
Signoravitch (2006) uses invariance arguments
to find powerful test statistics for multiple testing
when the underlying tests are multivariate.

\begin{figure}
\pspicture(0,0)(16,6)
\rput[bl](3,0){\includegraphics[angle=-90,width=4in]{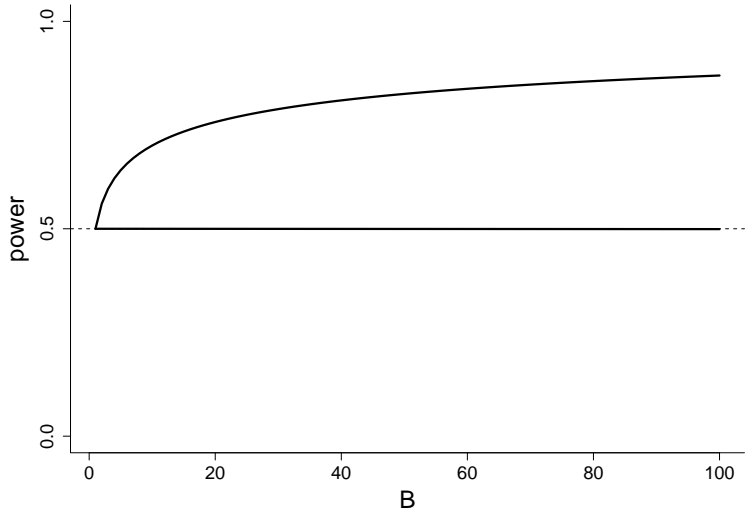}}
\endpspicture
\caption{Power gain/loss for weighting a single hypothesis.
In this example, an unweighted hypothesis has power 1/2.
The weights are $w_0 < 1 < w_1$ with $w_1/w_0 =B$.
The top line shows the power 
when the alternative is given the correct weight $w_1$.
The bottom line, which is nearly indistinguishable from 1/2,
shows the power 
when the alternative is given the incorrect weight $w_0$.
As $B$ increases, the power gain increases sharply while the power loss
remains nearly constant.}
\label{fig::single}
\end{figure}

\begin{figure}
\pspicture(0,0)(16,12)
\rput[bl](2,0){\includegraphics[width=5in,height=6in]{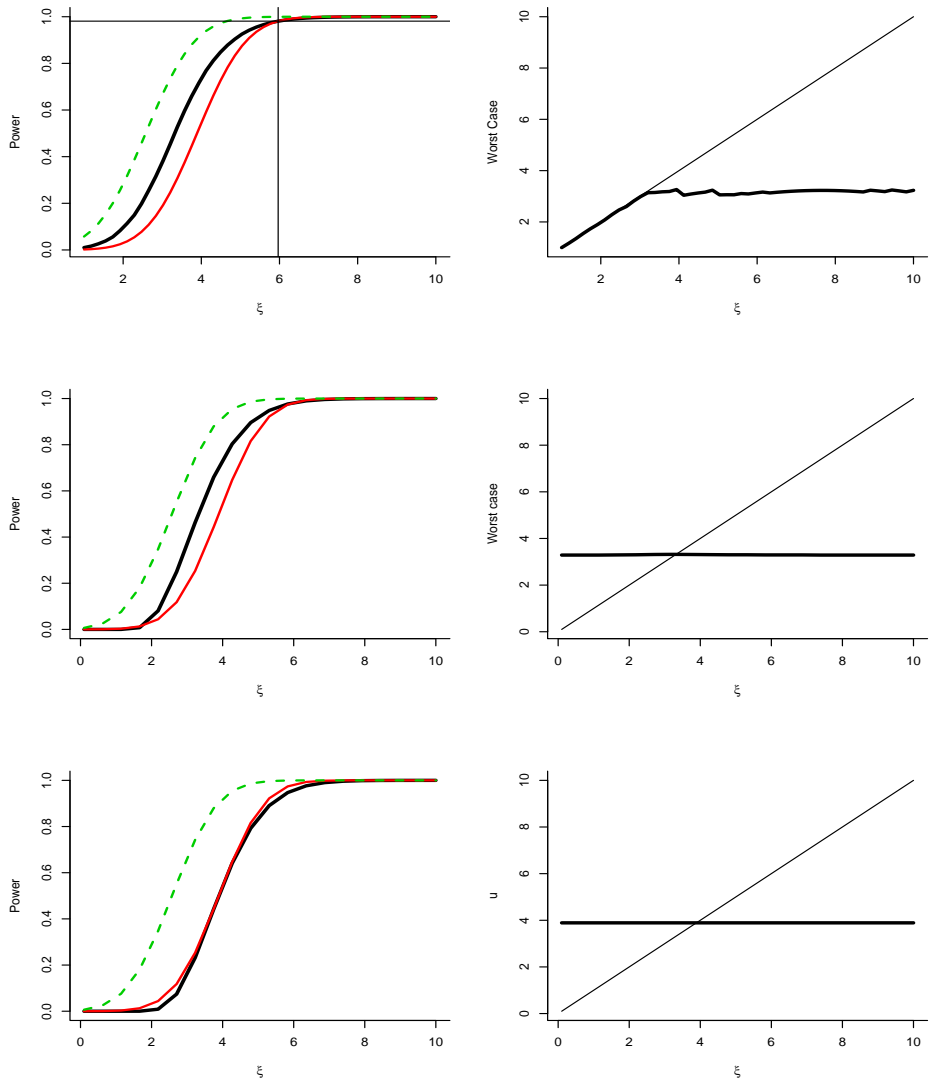}}
\endpspicture
\caption{Power as a function of the alternative mean $\xi$.
In these plots, $a=.01$, $m=1000$ and $\alpha=0.05$.
There are $(1-a)m$ nulls and $ma$ alternatives with mean $\xi$.
The left plots shows what happens when the weights are
incorrectly computed assuming that a fraction $\gamma$ of nulls
are actually alternatives with mean $u$.
In the top plot, we restrict $0 < u < \xi$.
In the second and third plot, no restriction is placed on $u$.
The top and middle plot have $\gamma = .1$ while
the third plot has $\gamma =1-a$ (all nulls misspecified as alternatives).
The dark solid line shows the lowest possible power assuming
the weights were estimated as poorly as possible.
The lighter solid line is the power of the unweighted (Bonferroni) method.
The dotted line is the power under the optimal weights.
The vertical line in the top plot is at $\xi_*$. The weighted method
beats unweighted for al $\xi < \xi_*$.
The right plot shows the least favorable $u$ as a function of $\xi$.
That is, mistaking $\gamma m$ nulls for alternatives with mean $u$
leads to the worst power.
Also shown is the line $u=\xi$.}
\label{fig::minimax}
\end{figure}

In this paper we show that the optimal weights form a one parameter
family.  We also show the power is very robust to misspecification of
the weights.  In particular, we show that (i) sparse weights 
(a few large weights and minimum weight close to 1)
lead to huge power gain for well specified weights, but minute power loss for
poorly specified weights; and (ii) in the non-sparse case, under weak
conditions, the worst case power loss for poorly specified weights is
typically better than the power using equal weights.  
In fact, the power is degraded at most by a factor of about
$\gamma/(1-a)$
where $a$ is the fraction of nonnulls and $\gamma$ is the fraction of nulls
that are mistaken for alternatives.
Figure \ref{fig::single} shows the sparse case.  The top line shows power
from correct weighting while the bottom line shows power from incorrect
weighting.  
We see that the power gains overwhelm the potential power loss.
Figure \ref{fig::minimax} shows the non-sparse case.  The
plots on the left show the power as a function of the alternative mean $\xi$.  The
dark solid line shows the lowest possible power assuming the weights
were estimated as poorly as possible.  The lighter solid line is the
power of the unweighted (Bonferroni) method.  
The dotted line shows the power under theoretically optimal weights.
The worst case weighted
power is typically close to or larger then the Bonferroni power except for large $\xi$
when they are both large.

We consider two methods for choosing the weights: (i) external
weights, where prior information (based on scientific knowledge or
prior data) singles out specific hypotheses and (ii) estimated weights
where the data are used to construct weights.  External weights are
prone to bias while estimated weights are prone to variability.  The
two robustness properties reduce concerns about bias and variance.

An example of external weighting is the following.  We have test
statistics $\{T_j:\ j=1, \ldots, m\}$ associated with spatial
locations $\{s_j:\ j=1, \ldots, m\}$ where $s_j \in [0,L]$, say.  
These could be association tests for markers on a genome.
The number of tests $m$ is large, on the order of 100,000 for example.
Each $T_j$ is used to test the null hypothesis that $\theta_j =
\E(T_j)=0$.  Prior data is in the form of a smooth stochastic process
$\{Z(s): s\in [0,L]\}$.  
This might be from a whole genome linkage scan.
At alternatives, the mean $\mu(s)=\E(Z(s))$ is a large 
positive value; however, due to correlation, at nulls close to
alternatives, $\mu(s)$ is also non-zero.  Peaks in the process $Z(s)$
provide approximate information about the location of alternatives.
We want to use the process $Z$ to generate reasonable weights for the
test statistics.

When external weights are not available, the optimal weights
can be estimated from the data.
One approach is to use data splitting (RVD) using a fraction of the data to
estimate the weights and the remainder to test.
For example, consider the two-stage genome-wide
association study (e.g., Thomas et al. 2005) for which a sample of $n$
subjects is split into two subsets. Using the first subset, we obtain
test statistics $\{T_j:\ j=1, \ldots, m\}$ associated with 
locations $\{s_j:\ j=1, \ldots, m\}$.
Typically only the second subset of data are used in the
final analysis. Building on the ideas of Skol et al. (2006),
we take the two-stage
study design further, exploring how the first set of data can be
utilized to formulate weights, and the full data set can be used for testing.

\section{Weighted Multiple Testing}

We are given hypotheses $H = (H_1, \ldots, H_m)$ and standardized test
statistics $T = (T_1, \ldots, T_m)$ where $T_j \sim N(\xi_j,1)$.
(The methods can be extended for nonnormal test statistics
but we do not consider that case here.)
For a two-sided hypothesis, $H_j =1$ if $\xi_j \neq 0$ and $H_j =0$
otherwise.  For the sake of parsimony, unless otherwise noted, 
results will be stated for a one-sided test where $H_j =1$ if
$\xi_j > 0$ although the results extend easily to the two-sided case.
Let $\theta = (\xi_1, \ldots, \xi_m)$ denote the vector of means.

The original data are often of the form
\begin{equation}
\label{eq:xmat}
\mathbb{X} = \left(\begin{array}{cccc}
X_{11} & X_{12} & \ldots & X_{1m}\\
X_{21} & X_{22} & \ldots & X_{2m}\\
\vdots & \vdots & \vdots & \vdots \\
X_{n1} & X_{n2} & \ldots & X_{nm}
\end{array}\right)
\end{equation}
where the $j^{\rm th}$ test statistic $T_j$ is based on the $j^{\rm
  th}$ column of $\mathbb{X}$.  Usually, $T_j$ is of the form $T_j =
\sqrt{n_j}\overline{X}_j/\sigma_j$ where $\overline{X}_j$ is
approximately (or exactly) $N(\gamma_j, \sigma_j2/n_j)$ and the
noncentrality parameter is $\xi_j = \sqrt{n}_j\gamma_j/\sigma_j$.

The p-values associated with the tests are $P = (P_1, \ldots, P_m)$
where $P_j = \overline{\Phi}(T_j)$, $\overline{\Phi} = 1- \Phi$ and
$\Phi$ denotes the standard Normal {\sc cdf}.  Let
$$
P_{(1)} \leq \cdots \leq P_{(m)}
$$ 
denote the sorted p-values and let
$$
T_{(1)} \geq \cdots \geq T_{(m)}
$$ 
denote the sorted test statistics.

A {\em rejection set} ${\cal R}$ is a subset of $\{1, \ldots, m\}$.
Say that ${\cal R}$ {\em controls familywise error at level $\alpha$}
if $\P({\cal R}\cap {\cal H}_0)\leq \alpha$ where ${\cal H}_0 = \{j:\
H_j =0\}$.  The {\em Bonferroni rejection} set is
\begin{equation}
{\cal R} = \{j:\ P_j < \alpha/m\} = \{j: T_j > z_{\alpha/m}\}
\end{equation}
where we use the notation 
$z_{\beta} = \overline{\Phi}^{-1}(\beta)$.

The weighted Bonferroni procedure of 
Genovese, Roeder and Wasserman (2005) 
is as follows.
Specify nonnegative weights $w=(w_1, \ldots, w_m)$ and
reject hypothesis $H_j$ if 
\begin{equation}
j\in {\cal R} = \left\{ j:\ \frac{P_j}{w_j} \leq \frac{\alpha}{m}\right\}.
\end{equation}
As long as $m^{-1}\sum_j w_j =1$,
the rejection set ${\cal R}$
controls familywise error at level $\alpha$.
For completeness, we provide the proof.
(All further proofs are in the appendix.)

\begin{lemma}
If $m^{-1}\sum_j w_j =1$,
then the rejection set ${\cal R}$
controls familywise error at level $\alpha$.
\end{lemma}

\noindent
{\sf Proof.}
The familywise error is
\begin{eqnarray*}
\P(({\cal R}\cap {\cal H}_0) > 0)
  &  =   & \P\left( P_j \le \frac{\alpha w_j }{m}\ \ {\rm for\ some\ }j \in {\cal H}_0\right)\\
  & \leq & \sum_{j\in {\cal H}_0} \P\left( P_j \le \frac{\alpha w_j }{m}\right)=
  \frac{\alpha}{m} \sum_{j\in {\cal H}_0} w_j \leq \alpha \overline{w} = \alpha .\ \ \ \blacksquare
\end{eqnarray*}

Genovese, Roeder and Wasserman (2005) also showed that
false discovery methods benefit by weighting.
Recall that the false discovery proportion (FDP) is
\begin{equation}
{\rm FDP} = \frac{{\rm number\  of\  false\ rejections}}{{\rm number\  of\  rejections}} =
\frac{| {\cal R} \cap {\cal H}_0 |}{| {\cal R}|}
\end{equation}
where the ratio is defined to be 0 if the denominator is 0.
The false discovery rate (FDR) is
${\rm FDR} = \E({\rm FDP})$.
Benjamini and Hochberg (1995) proved
${\rm FDR} \leq \alpha$ if
${\cal R} = \{ j:\ P_{(j)} \leq T\}$
where
$T = \max\{ j:\ P_{(j)} \leq j\alpha/m\}$.
Genovese, Roeder and Wasserman (2004) showed that ${\rm FDR} \leq
\alpha$ if the $P_j's$ are replaced by $Q_j = P_j/w_j$ as long as
$m^{-1}\sum_j w_j =1$ as before.
This paper will focus only on familywise error.
Similar results hold for FDR
and will be in a followup paper.

\section{Power and Optimality}

\subsection{Power}

Before weighting, that is using weight 1, the power of a single, one-sided alternative is
\begin{equation}
\pi(\xi_j,1) = \P(T_j > z_{\alpha/m}) = \overline{\Phi}(z_{\alpha/m} - \xi_j).
\end{equation}
The power\footnote{
For a two-sided alternative the power is
$$
\pi(\xi_j,w_j)= \overline{\Phi}
\left(\overline{\Phi}^{-1}\left(\frac{\alpha w_j}{2m}\right)-\xi_j\right) +
\overline{\Phi}
\left(\overline{\Phi}^{-1}\left(\frac{\alpha w_j}{2m}\right)+\xi_j\right).
$$}
in the weighted case is
\begin{equation}\label{eq::power-formula}
\pi(\xi_j,w_j) =
\P\left( P_j < \frac{\alpha w_j}{m}\right) 
= \P\left( T_j > \overline{\Phi}^{-1}\left(\frac{\alpha w_j}{m}\right)\right) =
 \overline{\Phi}
\left(\overline{\Phi}^{-1}\left( z_{\alpha w_j/m}\right)-\xi_j\right).
\end{equation}
Weighting increases the power when $w_j >1$ and decreases the power when $w_j < 1$.

Given $\theta = (\xi_1, \ldots, \xi_m)$ and
$w = (w_1, \ldots, w_m)$ we define the {\em average power} 
\begin{equation}
\frac{1}{m}\sum_{j=1}^m \pi(\xi_j,w_j) I(\xi_j > 0).
\end{equation}
More generally, if $\xi$ is drawn from a distribution $Q$
and $w = w(\xi)$ is a weight function
we define the average power 
\begin{equation}
\int  \pi(\xi,w(\xi))I(\xi>0) dQ(\xi).
\end{equation}
If we take $Q$ to be the empirical distribution of
$(\xi_1, \ldots, \xi_m)$ 
then this reduces to the previous expression.
In this case we require $w(\xi) \geq 0$ and $\int w(\xi)dQ(\xi) =1$.

\subsection{Optimality and Robustness}

In the following theorem 
we see that the set of optimal weight functions form a one parameter
family indexed by a constant $c$.

\begin{theorem}\label{thm::optimal}
Given
$\theta=(\xi_1, \ldots, \xi_m)$,
the optimal weight 
vector $w=(w_1, \ldots, w_m)$ that maximizes the average power
subject to $w_j \geq 0$ and $m^{-1}\sum_{j=1}^m w_j=1$
is $w=(\rho_c(\xi_1), \ldots, \rho_c(\xi_m))$
where
\begin{equation}
\label{eq:wq1}
\rho_c(\xi) = 
\left(\frac{m}{\alpha}\right)\overline{\Phi}\left(\frac{\xi}{2} + \frac{c}{\xi}\right)I(\xi > 0),
\end{equation}
and $c\equiv c(\theta)$ is defined by the condition
\begin{equation}
\label{eq:wbest}
\frac{1}{m} \sum_{j=1}^m \rho_c(\xi_j)  =1.
\end{equation}
\end{theorem}

\begin{figure}
\pspicture(0,0)(16,8)
\rput[bl](3,0){\includegraphics[width=4in,height=4in]{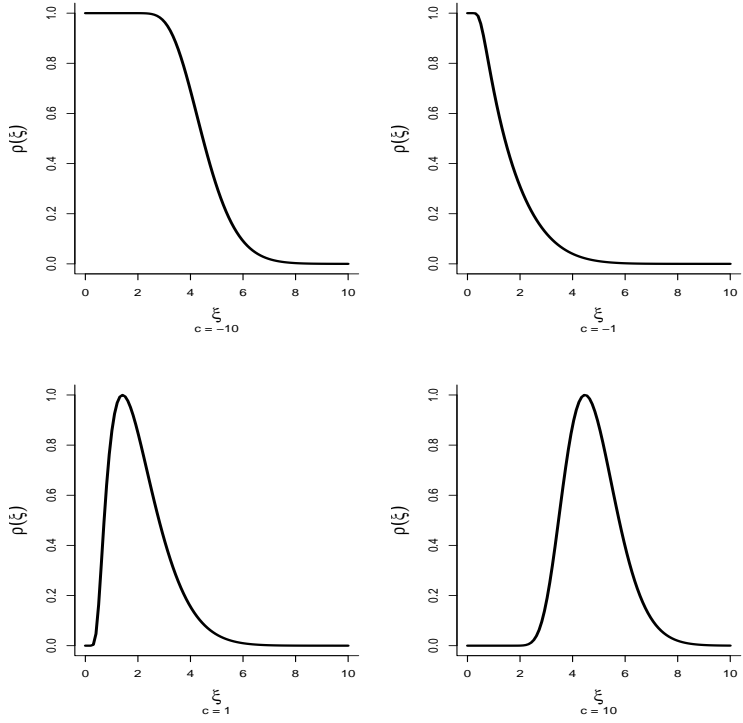}}
\endpspicture
\caption{Optimal weight function $\rho_c(\xi)$ for various $c$.
In each case $m=1000$ and $\alpha = 0.05$.
The functions are normalized to have maximum 1.}
\label{fig::family}
\end{figure}

The proof is in the appendix.
Some plots of the function $\rho_c(\xi)$ for various values of $c$ are 
shown in Figure \ref{fig::family}.
In these plots, the function is normalized to have maximum 1
for easier visualization.
The result generalizes to the case where the alternative means
are random variables with distribution $Q$ in which case
$c$ is defined by
\begin{equation}
\int \rho_c(\xi)dQ(\xi)  =1.
\end{equation}

{\bf Remark.}
Rejecting when $P_j/w_j \leq \alpha/m$ is the same as rejection when
$Z_j > \xi_j/2 + c/\xi_j$.
This is identical to the result
of Rubin, van der Laan and Dudoit (2005),
obtained independently.
The remainder of the paper,
which shows some good properties
of the weighted method, can thus also be
considered as providing support for
their method.
In particular, they noted in their
simulations then even poorly specified estimates of the cutoffs
$\xi_j/2 + c/\xi_j$ can still perform well.
This paper provides insight into why that is true.

\relax From (\ref{eq::power-formula}) and 
(\ref{eq:wq1}) we have immediately:

\begin{lemma}\label{lemma::power-at-opt}
The power at an alternative with mean $\xi$ under 
optimal weights is
$\overline{\Phi}\left(c/\xi - \xi/2\right)$.
The average power under optimal weights, which we call the
{\bf oracle power}, is
\begin{equation}
\frac{1}{m}\sum_{j=1}^m \overline{\Phi}\left(\frac{c}{\xi_j} - \frac{\xi_j}{2}\right)I(\xi_j > 0).
\end{equation}
\end{lemma}

The oracle power is not attainable since
the optimal weights depend on
$\theta=(\xi_1, \ldots,\xi_m)$ or,
equivalently, on $Q$.
In practice, the weights will either be chosen by
prior information or by estimating the $\xi$'s.
This raises the following question:
how sensitive is the power to correct specification of
the weights?
Now we show that the power is very robust to weight misspecification.

The weights themselves can be very sensitive
to changes in $\theta$.
Consider the following example.
Suppose that $\theta=(\xi_1, \ldots, \xi_m)$
where each $\xi$ is equal to either 0 or some fixed number $\xi$.
The empirical distribution of the $\xi_j$'s is thus
$Q=(1-a)\delta_0 + a \delta_\xi$
where $\delta$ denotes a point mass and $a$ is the fraction
of nonzero means.
The optimal weights are
$0$ for $\xi_j=0$ and 
$1/a$ for $\xi_j=\xi$.
Let
$\tilde{Q}=(1-a-\gamma)\delta_0 + \gamma\delta_u + a \delta_\xi$
where $u$ is a small positive number.
Since we have only moved the mass at 0 to $u$, and $u$
is small, we would hope that $w(\xi)$
will not change much.
But this is not the case.
Set
\begin{equation}
\xi    = A + \sqrt{A^2 - 2c},\ \ \ \ 
u     = B - \sqrt{B^2 - 2c}
\end{equation}
where
\begin{equation}
A     = \overline{\Phi}^{-1}\left(\frac{\alpha}{(m(\gamma K+a))}\right),\ \ \ \ 
B     = \overline{\Phi}^{-1}\left(\frac{K \alpha}{(m(\gamma K+a))}\right),
\end{equation}
yields weights $w_0$ and $w_1$ on $u$ and $\xi$ such that
$w_0/w_1=K$.
For example, take
$m= 1000$,  $\alpha = 0.05$,
$a     = .1$,
$\gamma   = .1$,
$K     = 1000$,
and $c     = .1$.
Then $u=.03$ and $\xi = 9.8$.
The optimal weight on $\xi$ under $Q$ is
10 but under $\tilde{Q}$ it is
$.00999$ and so is reduced by a factor of 1001.
More generally we have the following result
which shows that the weights are,
in a certain sense, a discontinuous function of $\theta$.

\begin{lemma}\label{lemma::discon}
Fix $\alpha$ and $m$.
For any $\delta>0$ and $\epsilon>0$
there exists
$Q=(1-a)\delta_0 + a\delta_\xi$ and
$\tilde{Q}=(1-a-\gamma)\delta_0  +\gamma \delta_u + a\delta_\xi$
such that
\begin{equation}
d(Q,\tilde{Q}) < \delta,\ \ \ {\rm and}\ \ \ 
\frac{\tilde\rho(\xi)}{\rho(\xi)} < \epsilon
\end{equation}
where
$a=\alpha/4$,
$d(Q,\tilde{Q}) = \sup_{\xi} | Q(-\infty, \xi], \tilde{Q}(-\infty, \xi]|$
is the Kolmogorov-Smnirnov distance,
$\rho$ is the optimal weight function for $Q$ and
$\tilde{\rho}$ is the optimal weight function for $\tilde{Q}$.
\end{lemma}

Fortunately, this problem is not serious
since it is possible to have high power even with poor weights.
In fact, the power of the weighted method has the following 
two robustness properties:

\vspace{1cm}

\noindent
{\bf Property I:}
Sparse weights (minimum weight close to 1) are highly robust.
If most weights are less than 1 and the minimum weight is close to 1 then
correct specification (large weights on alternatives) leads to large power gains but
incorrect specification (large weights on nulls) leads to little power loss.

\vspace{1cm}

\noindent
{\bf Property II:}
Worst case analysis.
Weighted hypothesis testing,
even with poorly chosen weights,
typically does as well or better than
Bonferroni except when the
the alternative means are large,
in which both have high power.

\vspace{1cm}

Let us now make the these statements precise.
Also, see Genovese, Roeder and Wasserman (2006) and
Roeder, Bacanu, Wasserman and Devlin (2006) for other results on the
effect of weight misspecification.

\vspace{1cm}

\noindent
{\bf Property I.}
Consider first the case where the weights take two distinct values
and the alternatives have a common mean $\xi$.
Let $\epsilon$ denote the 
fraction of hypotheses given the 
larger of the two values of the weights $B$.
Then, the weight vector $w$ is proportional to
$$
(\underbrace{B, \ldots, B}_{k\ {\rm terms}}, 
    \underbrace{1, \ldots, 1}_{m-k\ {\rm terms}})
$$
where $k = \epsilon m$ and $B > 1$
and hence the normalized weights are
\[
w = (\underbrace{w_1, \ldots, w_1}_{k\ {\rm terms}}, 
    \underbrace{w_0, \ldots, w_0}_{m-k\ {\rm terms}})
\]
where
$$
w_1 = \frac{B}{ \epsilon B + (1-\epsilon)},\ \ \ \ \ w_0 = \frac{1}{ \epsilon B + (1-\epsilon)}.
$$
We say that the weights are sparse if $\epsilon$ is small, that is,
if most weights are near 1.

\begin{figure}
\begin{center}
\includegraphics[height=3in,width=4in]{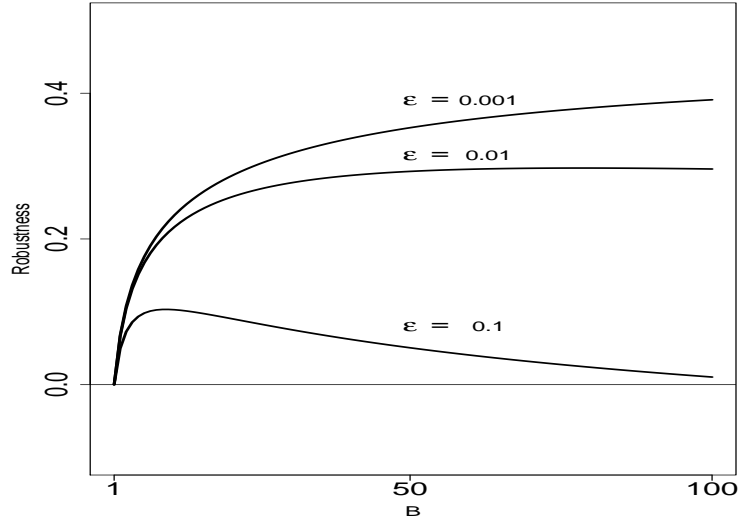}
\end{center}
\caption{Robustness function for $m=1000$.
In this example, $\xi = z_{\alpha/m}$ which has
power 1/2 without weighting.
The gain of correct weighting far outweighs the loss for
incorrect weighting as long as the fraction of large weights $\epsilon$
is small.}
\label{fig::powermult}
\end{figure}

Consider an alternative with mean $\xi$.
The power gain by correct weighting
is the power under weight $w_1$ minus the unweighted power 
$\pi(\xi,w_1)-\pi(\xi,1)$.
Similarly, the power loss for incorrect weighting is
$\pi(\xi,1)-\pi(\xi,w_0)$.
The gain minus the loss,
which we call the robustness function,  is
\begin{eqnarray}
R(B,\epsilon) & \equiv &
\biggl(\pi(\xi,w_1) -   \pi(\xi,1)\biggr) + \biggl(\pi(\xi,1) - \pi(\xi,w_0)\biggr)\\
&=&
\overline{\Phi}\left( z_{\alpha w_1/m} - \xi  \right) +
\overline{\Phi}\left( z_{\alpha w_0/m} - \xi  \right) -
2\overline{\Phi}\left( z_{\alpha /m} - \xi  \right).
\end{eqnarray}
The gain outweighs the loss if and only if
$R(B,\epsilon) >0$.
This is illustrated in
figures 
\ref{fig::single} and
\ref{fig::powermult}.

\begin{theorem}
Fix $B>1$.
Then,
$\lim_{\epsilon\to 0}R(B,\epsilon) > 0.$
Moreover, 
there exists $\epsilon^*(B)>0$ such that
$R(B,\epsilon)>0$ for all $\epsilon < \epsilon^*(B)$.
\end{theorem}

We can generalize this beyond the two-valued case as follows.
Let $w$ be any weight vector such that $m^{-1}\sum_jw_j =1$.
Now define the (worst case) robustness function
\begin{equation}
R(\xi)  \equiv 
\min_{\{j:\ w_j > 1, H_j=1\}}
\left\{ \pi(\xi,w_j) - \pi(\xi,1)\right\} - 
\max_{\{j:\ w_j < 1, H_j=1\}}
\left\{ \pi(\xi,1) - \pi(\xi,w_j)\right\} .
\end{equation}
We will see that $R(\xi) >0$ under weak conditions and that
the maximal robustness is obtained for $\xi$ near
the Bonferroni cutoff $z_{\alpha/m}$.

\begin{theorem}\label{thm::sparse-is-good}
A necessary and sufficient condition for $R(\xi) >0$ is
\begin{equation}\label{eq::cond-for-robust}
R_{b,B}(\xi) \equiv {\Phi}\Biggl( z_{\alpha B/m}-\xi\Biggr) +
{\Phi}\Biggl( z_{\alpha b/m}-\xi\Biggr)-
2{\Phi}\Biggl( z_{\alpha/m}-\xi \Biggr)\leq 0
\end{equation}
where
$B= \min\{w_j:\ w_j > 1\}$,
$b = \min\{w_j\}$.
Moreover,
\begin{equation}
R_{b,B}(\xi) = -\Delta(\xi) + O(1-b)
\end{equation}
where
\begin{equation}
\Delta(\xi) =  \left({\Phi}\Biggl( z_{\alpha/m}-\xi \Biggr)-
                       {\Phi}\Biggl( z_{\alpha B/m}-\xi\Biggr)\right)  > 0
\end{equation}
and, as $b\to 1$,
$\mu(\{\xi:\ R(\xi) < 0\})\to 0$ and
$\inf_\xi R(\xi) \to 0$.
\end{theorem}

\begin{figure}
\begin{center}
\includegraphics[angle=-90,width=7in]{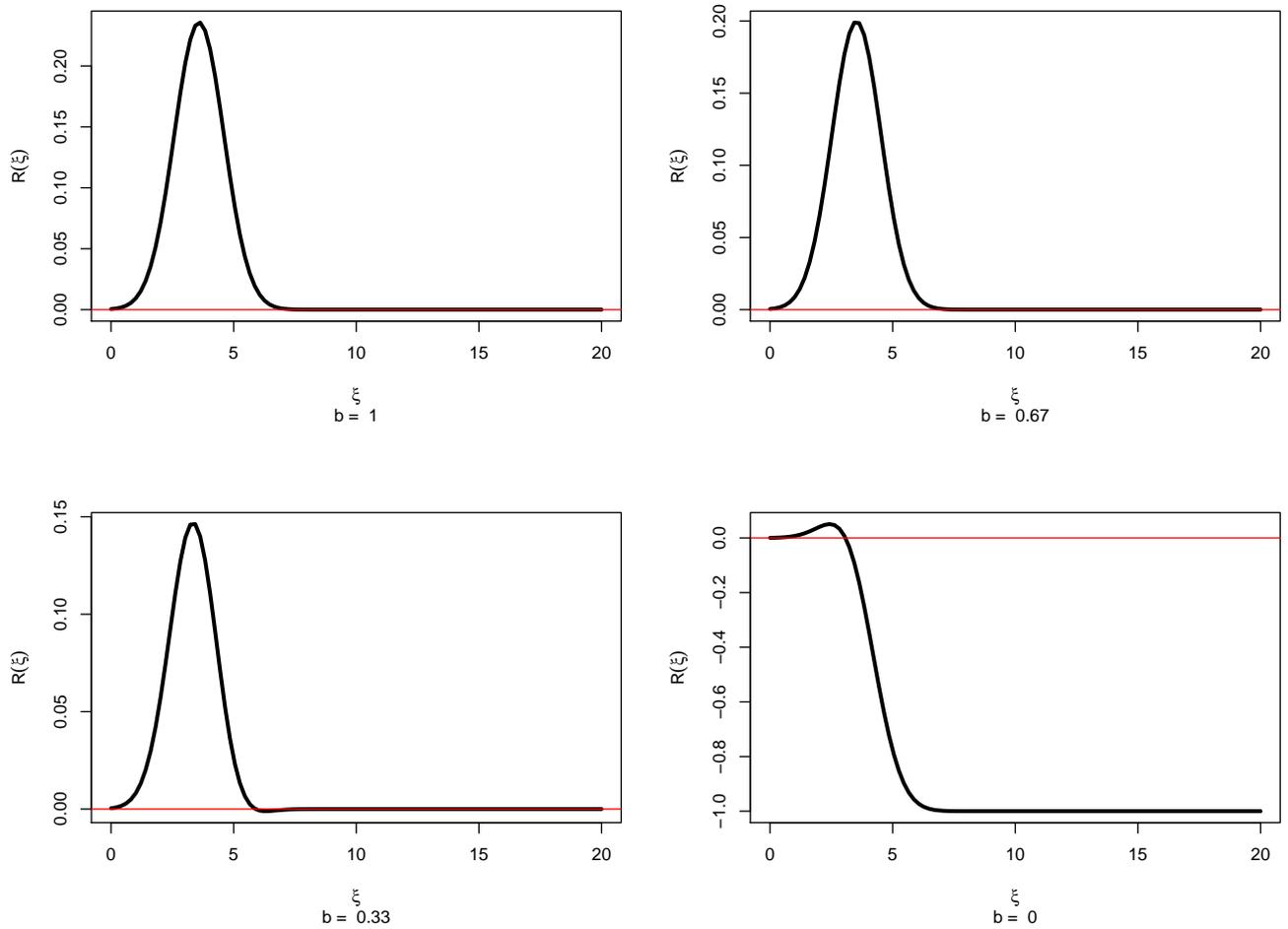}
\end{center}
\caption{The robustness function $R(\xi)$ for several values
of $b=\min_j w_j$. In each case, $m=1000$, $\alpha=0.05$, $B=10$.
Whenever $R(\xi) >0$, power gain outweighs power loss.
When $b$ is near 1, $R(\xi) >0$ for most $\xi$.
Even when $b=0$ there is a safe zone including $\xi=0$ as long as $B\geq 2$.}
\label{fig::robust}
\end{figure}

The theorem is illustrated in Figure 
\ref{fig::robust}.
We see that there is overwhelming robustness as long as the minimum weight
is near 1.
Even in the extreme case $b=0$,
there is still a safe zone, an interval of values of $\xi$ over which $R(\xi)>0$.

\begin{lemma}\label{lemma::safe}
Suppose that $B\geq 2$.
Then there exists $\xi_* > 0$ such that
$R_{B,b}(\xi) >0$ for all $0\leq \xi \leq \xi_*$ and all $b$.
An upper bound on $\xi_*$ is
$z_{\alpha/m} - 1/(z_{\alpha/m} - z_{B \alpha/m})$.
\end{lemma}

\vspace{1cm}

\noindent
{\bf Property II.}
Even if the weights are not sparse, 
the power of the weighted test cannot be too bad as we now show.
To begin, assume that each mean is either equal to $0$ or
$\xi$ for some fixed $\xi>0$.  
Thus, the empirical distribution is
\begin{equation}
Q = (1-a)\delta_0 + a \delta_\xi
\end{equation}
where $\delta$ denotes a point mass and $a$ is the
fraction of nonzero $\xi_j$'s.  The optimal weights are 
$1/a$ for hypotheses whose mean is $\xi$.
To study
the effect of misspecification error, consider the 
case where 
$b=\gamma m$ nulls are mistaken for alternatives with mean $u>0$.
This corresponds to misspecifying $Q$ to be
\begin{equation}
 \tilde{Q} = (1-a-\gamma)\delta_0 + \gamma \delta_u + a \delta_\xi.
\end{equation}
We will study the effect of varying $u$ so
let $\pi(u)$ denote the power at the true alternative $\xi$ as a function of $u$.
Also, let $\pi_{\rm Bonf}$ denote the power
using equal weights (Bonferroni).
Note that changing
$Q = (1-a)\delta_0 + a \delta_\xi$ to
$Q = (1-a)\delta_0 + a \delta_{\xi'}$  for $\xi' \neq \xi$
does not change the weights.

As the weights are a function of $c$, we first need to find 
$c$ as a function of $u$.
The normalization condition
(\ref{eq:wbest}) reduces to
\begin{equation}\label{eq::wbest2}
\gamma\overline{\Phi}\left(\frac{u}{2} + \frac{c}{u}\right) + 
a \overline{\Phi}\left(\frac{\xi}{2} + \frac{c}{\xi}\right) =\frac{\alpha}{m}
\end{equation}
which implicitly defines the function $c(u)$.
First we consider what happens when $u$ is restricted to be less than $\xi$.

\begin{theorem}\label{thm::restrict-u}
Assume that $\alpha/m \leq \gamma + a \leq 1$.
Let $Q = (1-a)\delta_0 + a \delta_\xi$ and
$\tilde{Q} = (1-a-\gamma)\delta_0 + \gamma \delta_u + a \delta_\xi$
with $0\leq u \leq \xi$.
Let $C(\xi) = \sup_{0\leq u \leq \xi}c(u)$ and define
$\xi_0 = z_{\alpha/(m(\gamma +a))}$,
\begin{enumerate}
\item For $\xi \leq \xi_0$,
\begin{equation}
C(\xi) = \xi \xi_0- \xi^2/2.
\end{equation}
For $\xi > \xi_0$, $C(\xi)$ is the solution to
\begin{equation}
\gamma \overline{\Phi}(\sqrt{2c}) + a\overline{\Phi}\left(\frac{c}{\xi} + \frac{\xi}{2}\right) = 
\frac{\alpha}{m}.
\end{equation}
In this case,
$C(\xi) = z^2_{\alpha/(m\gamma)}/2 + O(a)$.
\item Let 
\begin{equation}
\xi_* = 
z_{\alpha/m} + 
\sqrt{z_{\alpha/m}^2 - z_q^2},\ \ \ {\rm where}\ q = \frac{\alpha(1-a)}{m\gamma}.
\end{equation}
For $\xi < \xi_*$,
\begin{equation}\label{eq::thebonfclaim}
\inf_{0 < u < \xi}\pi(u) \geq \pi_{\rm Bonf}.
\end{equation}
For $\xi \geq \xi_*$ we have
\begin{eqnarray}
\label{eq::thenextclaim}
\inf_{0 < u < \xi}\pi(u) &\geq &
\overline{\Phi}\left( \frac{ z^2_{\alpha/(m\gamma)} - \xi_*^2}{2\xi_*}\right)- O(a)\\
&\approx & 1-\overline{\Phi}\left(\sqrt{2\log\frac{1-a}{\gamma}}\right)-O(a)\\
& \geq & 1- \frac{\gamma}{1-a} - O(a).
\end{eqnarray}
\end{enumerate}
\end{theorem}

The factor
$\overline{\Phi}\left(\sqrt{2\log\frac{1-a}{\gamma}}\right)\approx \gamma/(1-a)$
is the worst case power deficit due to misspecification.

Now we drop the assumption that $u\leq \xi$.

\begin{theorem}\label{thm::drop-u}
Let $Q=(1-a)\delta_0 + a\delta_\xi$ and
let
$Q_u \equiv (1-a-\gamma)\delta_0 + \gamma\delta_u +a\delta_\xi$.
Let $\pi_u$ denote the power at $\xi$ using the weights
computed under $Q_u$.
\begin{enumerate}
\item
The least favorable $u$ is
\begin{equation}
u_* \equiv {\rm argmin}_{u\geq 0} \pi_u = \sqrt{2c_*} = z_{\alpha/(m\gamma)} + O(a)
\end{equation}
where
$c_*$ solves
\begin{equation}
\gamma\overline{\Phi}(\sqrt{2 c_*}) + a \overline{\Phi}\left( \frac{\xi}{2}+\frac{c_*}{\xi}\right)=
\frac{\alpha}{m}
\end{equation}
and
$c_* = z^2_{\alpha/(m\gamma)}/2 + O(a)$.
\item The minimal power is
\begin{equation}
\inf_u \pi_u = \overline{\Phi}\left(\frac{c_*}{\xi} - \frac{\xi}{2}\right) =
\overline{\Phi}\left(\frac{z_{\alpha/(m\gamma)}^2-\xi^2}{2\xi}\right) + O(a).
\end{equation}
\item A sufficient condition for 
$\inf_u \pi_u$ to be larger than the power of the Bonferroni method is
\begin{equation}
\xi \geq z_{\alpha/m} + \sqrt{z^2_{\alpha/m} - z^2_{\alpha/(m\gamma)}} + O(a).
\end{equation} 
\end{enumerate}
\end{theorem}

\section{Choosing External Weights}

In choosing external weights, we will focus here on the two-valued case.
Thus, 
\begin{equation}
w = (\underbrace{w_1, \ldots, w_1}_{k\ {\rm terms}}, 
    \underbrace{w_0, \ldots, w_0}_{m-k\ {\rm terms}})
\end{equation}
where $k=\epsilon m$,
$w_1 = B/( \epsilon B + (1-\epsilon))$ and
$w_0 = 1/( \epsilon B + (1-\epsilon))$.
In practice, we would typically have a fixed fraction of hypotheses $\epsilon$
that we want to give more weight to.
The question is how to choose $B$.
We will focus on choosing $B$ to produce weights with good properties
at interesting values of $\xi$.
Now large values of $\xi$ already have high power.
Very small values of $\xi$ have extremely low power and benefit little by weighting.
This leads us to focus on 
constructing weights that are
useful for a {\em marginal effect}, defined as the alternative $\xi_m$
that has power 1/2 when given weight 1.
Thus, the marginal effect is $\xi_m = z_{\alpha/m}$.
In the rest of this section then
we assume that all nonzero $\xi_j$'s are equal to $\xi_m$.
Of course, the validity of the procedure does not depend on this
assumption being true.

\begin{figure}
\begin{center}
\includegraphics[angle=-90,width=5in]{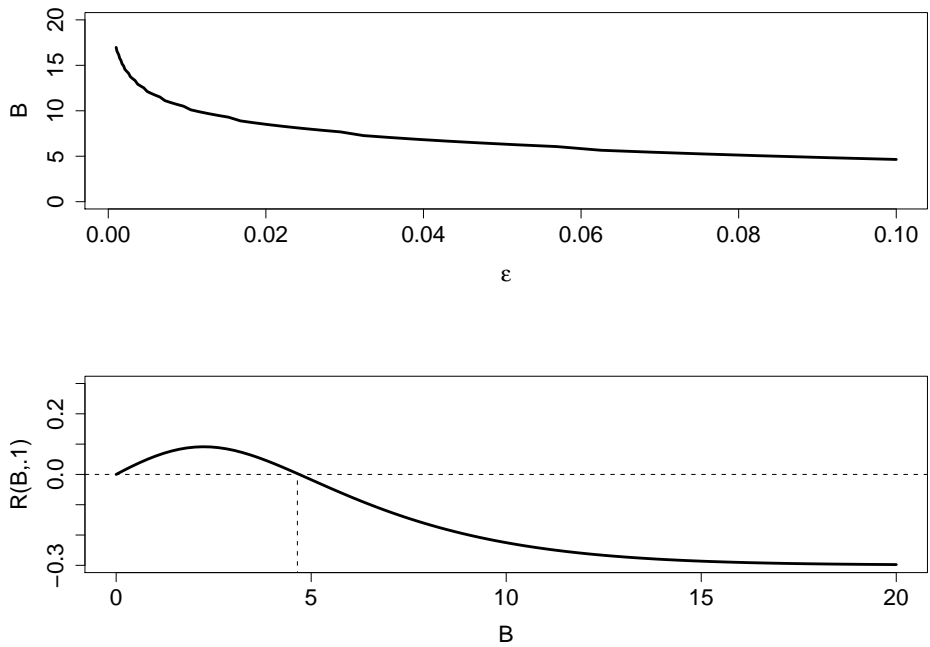}
\end{center}
\caption{Top plot: $B_0(\epsilon)$ versus $\epsilon$.
Bottom plot shows $R(B,.1)$
versus $B$. The turnaround point $B^*(\epsilon)$ is shown
with a vertical dotted line.}
\label{fig::turnaround}
\end{figure}

Fix $0 < \epsilon < 1$ and vary $B$.
As we increase $B$, we will eventually reach a point 
$B_0(\epsilon)$ where 
$R(B,\epsilon) < 0$ which we call
turnaround point.
Formally,
\begin{equation}
B_0(\epsilon) = \sup\biggl\{ B:\ R(B,\epsilon) > 0 \biggl\}.
\end{equation}
The top panel in Figure \ref{fig::turnaround}
shows 
$B_0(\epsilon)$ versus $\epsilon$
which shows that for small $\epsilon$ we can choose $B$ large
without loss of power.
The bottom panel shows
$R(B,\epsilon)$ for $\epsilon = 0.1$.
We suggest using 
$B=B_*(\epsilon)$, the value of $B$ that maximizes $R(B,\epsilon)$.

\begin{theorem}
Fix $0 < \epsilon < 1$.
As a function of $B$,
$R(B,\epsilon)$ is unimodal and satisfies
$R(1,\epsilon) =1$,
$R'(1,\epsilon) >0$ and
$R(\infty,\epsilon)  < 0$.
Hence,
$B_0(\epsilon)$ exists and is unique.
Also, $R(B,\epsilon)$ has a unique maximum at some point $B^*(\epsilon)$ and
$R(B^*(\epsilon),\epsilon)>0$.
\end{theorem}

When $\epsilon$ is very small, we can essentially choose $B$ as large as we like,
For example, suppose we want to
increase the chance of rejecting one particular hypothesis
so that $\epsilon = 1/m$.
Then,
\[
w_1 = \frac{mB}{B+m -1}\approx B,\ \ \ w_0 = \frac{1}{B+m -1}\approx 1
\]
and
\[
\lim_{m\to\infty}\lim_{B\to\infty} \pi(\xi_j,w_1) =1, \ \ \ 
{\rm while}\ \ \ \lim_{m\to\infty}\lim_{B\to\infty} \pi(\xi_j,w_0) =\frac{1}{2}.
\]
See Figure \ref{fig::single}.

The next results show that binary weighting schemes 
are optimal in a certain sense.
Suppose we want to have at least a fraction $\epsilon$
with high power $1-\beta$ and otherwise we want to
maximize the minimum power.

\begin{theorem}
\label{thm::method1}
Consider the following optimization problem:
Given $0 < \epsilon < 1$ and
$0 < \beta < 1/2$, find  a vector
$w = (w_1, \ldots, w_m)$ that maximizes
$$
\min_j \pi(\xi_m,w_j)
$$ 
subject to
$$
\overline{w} =1,\ \ \ \ \mbox{ and }\ 
\frac{\# \{ j:\ \pi(w_j,\xi_m)\geq 1-\beta\} }{m} \geq \epsilon.
$$
The solution is given by
\begin{equation}
w = (\underbrace{w_1, \ldots, w_1}_{k\ {\rm terms}}, 
    \underbrace{w_0, \ldots, w_0}_{m-k\ {\rm terms}})
\end{equation}
where
$w_1 = B/( \epsilon B + (1-\epsilon))$,
$w_0 = 1/( \epsilon B + (1-\epsilon))$,
$k = \epsilon m$,
$B =  cm(1-\epsilon)/(\alpha - \epsilon c m)$
and
$c = \overline{\Phi}\left( z_{\alpha/m} + z_{1-\beta}\right)$.
\end{theorem}

If our goal is to maximize the number of alternatives with
high power while maintaining a minimum power loss,
the solution is given as follows.

\begin{theorem}
\label{thm::method2}
Consider the following optimization problem:
Given $0 < \beta < 1/2$, find  a vector
$w = (w_1, \ldots, w_m)$ that maximizes
\begin{equation}
\# \{ j:\ \pi(w_j,\xi_m)\geq 1-\beta\}
\end{equation}
subject to
\begin{equation}
\overline{w} =1,\ \ \ {\rm and}\ 
\min_j \pi(w_j,\xi_m) \geq \delta .
\end{equation}
The solution is
\begin{equation}
w = (\underbrace{w_1, \ldots, w_1}_{k\ {\rm terms}}, 
    \underbrace{w_0, \ldots, w_0}_{m-k\ {\rm terms}})
\end{equation}
where
\begin{equation}
w_1 =\frac{m}{\alpha}\overline{\Phi}\left(z_{\alpha/m} + z_{1-\beta}\right),\ \ \ 
w_0 = \frac{m}{\alpha}\overline{\Phi}\left(z_{\alpha/m} + z_{\delta}\right),\ \ \ 
\epsilon = \frac{1-w_0}{w_1-w_0}
\end{equation}
and $k = m \epsilon$.
\end{theorem}

A special case that falls under this Theorem permits the minimum power
to be 0.  In this case $w_0=0$ and $\epsilon=1/w_1$.

\section{Estimated Weights}

In this section we explain how to use the data to estimate the weights.
There are two issues: we must ensure that
the error is still controlled and avoid
incurring large losses of power due to
replacing $\theta$ with an estimator $\hat\theta$.

\subsection{Validity With Estimated Weights}

{\bf Data Splitting.}  The approach, taken by RVD, for ensuring that
the error control is preserved relies on data splitting.  This
approach relies on normalized test statistics $T^{(l)},T^{(2)}$ based
on a partition of the data into subsets
$\mathbb{X}^{(1)},\mathbb{X}^{(2)}$ which include fractions $b$ and
$(1-b)$ of $\mathbb{X}$, respectively.  Note that $T_j =b^{1/2}\ta +
(1-b)^{1/2}\tb, j=1,\ldots,m$.  The training data $\mathbb{X}^{(1)}$
is used to estimate the noncentrality parameter of the standardized
statistic $T_j^{(1)}$, where $E[T_j^{(1)}]= \sqrt{b}\,\xi_j \equiv
\xi_j^{(1)}$.  Testing is conducted using the remaining fraction of
the data $\mathbb{X}^{(2)}$.  Consequently $\hat{\xi}_j^{(1)}$ must be
rescaled by $r_s = ({1-b})/{b})^{1/2}$ to estimate the
noncentrality parameter of the standardized statistic $T_j^{(2)}$,
i.e., $\hat{\xi}_j=r_s\,\hat{\xi}_j^{(1)}$.  The estimated weights are
$\hat{w}_j(T^{(1)}) =\rho_c(\hat{\xi}_j)$.  Because of the
independence between the two portions of the data, familywise error is
controlled at the nominal level.

\begin{lemma}
The procedure that rejects when
$P(T_j^{(2)}) < w(T^{(1)})\,\alpha/m$ controls the familywise
error at level $\alpha$.
\end{lemma}

\noindent
{\bf Recovering Power.}  As noted by Skol et al. (2005), data
splitting incurs a loss of power because the p-values are computed
using only a fraction the data.  To recover this lost power, we need
to use all the data to compute the p-values.  When using this approach
$\hat{\xi}_j^{(1)}$ must be rescaled by $r_f=b^{-1/2}$ to estimate the
noncentrality parameter of the standardized statistic $T_j$, i.e.,
$\hat{\xi}_j=r_f\,\hat{\xi}_j^{(1)}$. As in the data splitting
procedure, the estimated weights $\hat{w}_j(T^{(1)})
=\rho_c(\hat{\xi}_j)$ depend only on $\mathbb{X}^{(1)}$.  To preserve
error control we proceed as follows.

\begin{theorem}\label{eq::newcontrol}
Assume $T_{j}^{(k)}\sim N(0,1)$ independently for $k=1,2$.
Suppose that weight $w(T^{(1)})$ depends only on $\mathbb{X}^{(1)}$
but the p-value $P(T_j)$
is allowed to depend on the full data $\mathbb{X}$.
Define $c(T^{(1)})$ to solve
\begin{equation}
\frac{1}{m}\sum_{j=1}^m
\left(\overline{\Phi}
 \left(\frac{ \frac{\hat\xi_j}{2} + \frac{c(T^{(1)})}{\hat\xi_j} 
- \sqrt{b}\,T_j^{(1)} }  {\sqrt{1-b}}\right)\right) = \frac{\alpha}{m}.
\end{equation}
Then the procedure that rejects when
$P(T_j)  < w_j(T^{(1)})\alpha/m$, where
\begin{equation}
w_j(T^{(1)}) = \frac{m}{\alpha}\overline{\Phi}\left(\frac{\hat\xi_j}{2} 
+ \frac{c(T^{(1)})}{\hat{\xi_j}}\right)
\end{equation}
controls the familywise error at level $\alpha$.
\end{theorem}

\subsection{Simulations}

We simulate a study with $m=1000$ tests,
yielding data of the form given in (\ref{eq:xmat}).  A test
of the hypothesis $H_0:\xi_j\neq 0$ is performed for each $j$ using
$T_j$, which we assume is (approximately) normally distributed, or
equivalently $T_j^2 \sim \chi^2_1$.  In our simulations we generate 50
of the 1000 tests under the alternative hypothesis with shift
parameter $\xi_j =$ 2, 3, 4 or 5.  We compare the power for various
levels of a threshold parameter $\lambda \in (0,.5,1,1.5,2,2.5)$.  We
use a fraction $b=0.5$ of the data to construct the weights and we compare four
methods for estimating the noncentrality parameter:
\begin{enumerate}
\item The normalized statistic $\hat\xi_j^{(1)} = T^{(1)}_j$.
\item Hard thresholding:
\begin{equation}
\hat\xi_j^{(1)} =  T^{(1)}_j I(|T^{(1)}_j|> \lambda).
\end{equation}
\item Soft thresholding:
\begin{equation}
\hat\xi_j^{(1)} = {\rm sign}(T^{(1)}_j)(| T^{(1)}_j - \lambda)_+.
\end{equation}
\item The James-Stein estimator 
\begin{equation}
\hat\xi_j^{(1)} = \left(1 - \frac{m-2}{\sum_i  \left({T^{(1)}_i}\right)^2}\right)_+ T^{(1)}_j.
\end{equation}
To compute the weights we 
rescale $\hat\xi_j^{(1)}$ by $r_s$ or $r_f$ as appropriate to the followup testing
strategy.
\end{enumerate}

\begin{figure}
\begin{center}
\includegraphics[angle=-90,width=5in]{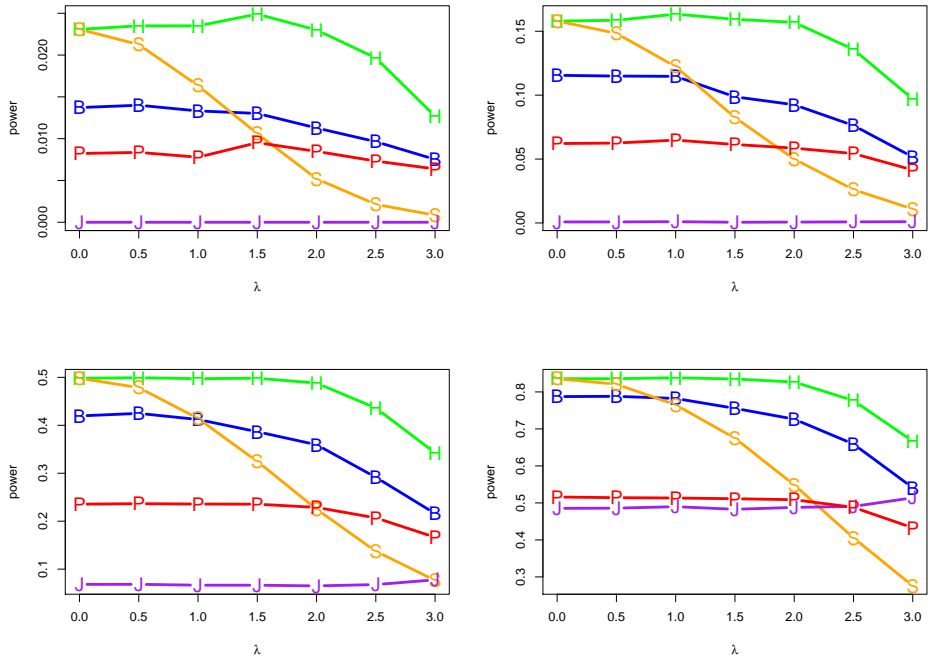}
\end{center}
\caption{Power of weighted tests.  From top left clockwise:
$\xi = 2,3,4,5$.
Methods compared use
weights based on hard thresholding (H), soft thresholding (S),
binary weights (B), and James-Stein (J).}
\label{fig::power}
\end{figure}

Power results are displayed in Fig. \ref{fig::power}.  We first
consider the power of the RVD procedure which uses the data splitting
strategy and $\lambda = 0$ (Fig. \ref{fig::power}, labeled ``P'' at
the origin).  Although RVD suggest using $\lambda = 0$, we also
examine the power of this procedure for a range of values of
$\lambda$.  This extended RVD procedure is applying hard-thresholding
to estimate $\theta$.  Next we consider the power of four testing
strategies that use the full data $T_{j}$ for testing rather than data
splitting.  The first approach (B) uses binary weights equal to $m/M$
where $M=\sum_i I\{|T_{j}^{(1)}|>\lambda\}$.  In this setting, when
$\lambda=0$ the method reduces to the simple one-stage Bonferroni
approach. For $\lambda>0$ it is the method of Skol et al. (2006).  The
remaining three approaches rely on weights estimated using
hard-thresholding (H), soft-thresholding (S), or James-Stein (J).  For
$\lambda=0$ methods H and S reduce to the normed sample mean which is
the RVD approach adapted to incorporate the full data in the p-value.
Clearly, this adaptation of the RVD method leads to a valuable
increase in power. For $\lambda>0$ this is method imposes a hard
threshold shrinkage effect on the parameter estimates.  Notice that
for any fixed value of $\lambda$, method H gives the best power.  In
particular, the difference in power between methods B and H
illustrates the advantage of using variable weights estimated from a
fraction of the data.  Method H and to a lesser extent method B are
nearly invariant to $\lambda$ for moderate values of the threshold
parameter.  In contrast, method S, relying on soft-thresholding,
experiences a sharp decline in power as $\lambda$ increases. Finally,
the James-Stein approach clearly fails in this setting, presumably
because most tests follow the null hypothesis and hence the true
signals are shrunk toward 0 which diminishes the power of the
procedure.

For each condition investigated the
tests had size less than 0.05
as expected from the theory.
The James-Stein method was most conservative.

\begin{figure}
\begin{center}
\includegraphics[angle=-90,width=5in]{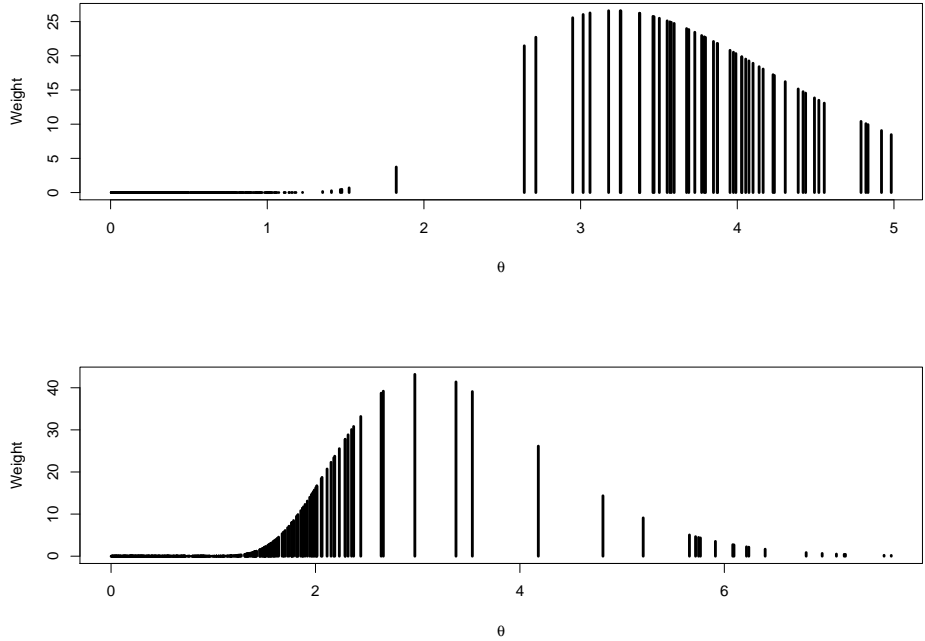}
\end{center}
\caption{Distribution of weights for two sets of data.}
\label{fig::wghtdist}
\end{figure}

\relax From this experiment it appears that shrinkage only enhances
power when the signal is very weak.  A more careful analysis reveals
that the effect of shrinkage for stronger signals is more subtle. As
$\xi\to 0$, $\rho_c(\xi)\to 0$.  Figure \ref{fig::wghtdist} shows
$\rho_c(\xi)$ is close to zero for a broad range of values.
Consequently, for $\lambda \leq 1.5$, the weight function performs
almost the same role as the threshold parameter.  Using
hard-thresholding for $\lambda < 1.5$ is essentially equivalent to
using using no threshold because a moderate level of shrinkage is
automatically imposed by the weight function.  Figure
\ref{fig::wghtdist} also illustrates how the optimal weights vary with
the signal strength (top panel has greater signal than bottom panel).
Both panels indicate that larger weights are placed in the midrange of
signal strength.  Essentially no weight is wasted on tests with small
signals ($\xi < 1.5$) because these tests are not likely to yield
significant results.  The bottom panel shows that large weights are
also not wasted on signals so strong that the tests can easily be
rejected even without up-weighting ($\xi >6$).  The top panel places
its largest weights between 2.5 and 4. The bottom panel has fewer
signals in this range and hence stronger weights can be applied to
signals between 2 and 2.5.  Both panels indicate near 0 weights would
be applied to tests with signals near 0.

\section{Discussion}

An interesting connection can be made between weights based on
threshold-estimators and two-stage experimental designs that perform
only a subset of the tests in stage two, based on the results obtained
from stage one.  The simplest example of this type of two-stage
testing is the two-stage Bonferroni
procedure, for which the training data $\mathbb{X}^{(1)}$ is used to
determine the $M$ elements in $\Lambda=\{j:|T_{j}^{(1)}|>\lambda\}$;
$\mathbb{X}^{(2)}$ is only measured for these columns.  A Bonferroni
correction with $\alpha/(2M)$ controls FWER at level $\alpha$ for
two-sided testing in this setting.  In essence this approach is a
weighted test with weights equal to $m/M$ for the elements in
$\Lambda$ and zero else where.

While the classic two-stage approach uses $\mathbb{X}^{(1)}$ for training,
and $\mathbb{X}^{(2)}$ for testing,
an alternative is to use the training data to determine the weights
and then use all of the data to conduct the tests.  This strategy was
recently investigated by Skol et al. (2006), using constant weights.
These authors use the training data to determine $\Lambda$ and then
apply weights equal to $m/M$ to the $M$ tests determined in stage one.
This full data approach proved to be
considerably more powerful than the two-stage Bonferroni
approach in simulations.

For hard and soft-thresholding, $\hat \xi_j=0$ for any
$|T_{j}^{(1)}|<\lambda$.  \relax From (\ref{eq:wq1}) it follows that
the weights for any test with $\hat \xi_j=0$ are 0 and the rejection
region is $Z_0=\infty$. Hence, a procedure using $w_j=0$ for columns with $ \hat
\xi_j=0$ is equivalent to a truncation procedure that tests only those
columns in $\Lambda$.  In practice, 
$\lambda$ can be chosen to optimize power or to constrain the experimental budget.  It
is worth noting that in some experimental settings, such as those
described by Skol et al., this experimental design can lead to
considerable savings of effort and resources.  Our results suggest
that this savings can be gleaned without loosing measurable power.

The same ideas used here can be applied to
other testing methods to improve power.
In particular, weights can be added to
the FDR method, Holm's stepdown test,
and the Donoho-Jin (2004) method.
Weighting ideas can also be used for
confidence intervals.
We plan to present the details for the other methods
in a followup paper.
Another item to be addressed in future work
is the connection with Bayesian methods.

As we noted, using weights is equivalent to using
a separate rejection cutoff for each statistic.
The methods of Storey (2005) and Signoravich (2006)
find optimal cutoffs when the cutoffs are constrained.
There is undoubtedly a bias-variance tradeoff.
These constrained methods
can estimate optimal cutoffs well (low variance)
but they will not achieve the oracle power obtained here
since they are by design biased away from these separate cutoffs.
Future work should be directed at comparing these approaches
and developing methods that lie in between these extremes.

\section{Appendix}

{\sf Proof of Theorem \ref{thm::optimal}.}
Let $A$ denote the set of hypotheses with $\xi_j > 0$.  Power is optimized
if $w_j = 0$ for $j \notin A$. The average power is
\[
\frac 1 m \sum_{j\in A} 
\phib\left(\phibi\left(\frac{\alpha w_j}m \right)-\xi_j\right).\]
with constraint 
\[ \sum_{j\in A}w_j = m.\]
Choose $\underline w$ to maximize 
\begin{eqnarray*}
\pi=\frac 1m\sum_{j \in A} \phib\left(\phibi\left(\frac{\alpha w_j}m\right) - \xi_j\right) 
- \lambda\left({m-\sum_{j\in A}w_i}\right)\\
\end{eqnarray*}
by setting the derivative to zero
\begin{eqnarray*}
\frac \partial{\partial w_i}\pi =   - \lambda  &+&
    \frac{\phi\left(\phibi\left(\frac{\alpha w_j}m\right) - \xi_j\right)}
         {\phi\left(\phibi\left(\frac{\alpha w_j}m\right)\right)} \frac{\alpha}m  = 0 \\
\frac{m\lambda}{\alpha } &=& \frac{\phi\left(\phibi\left(\frac{\alpha w_j}m\right) - \xi_j\right)}
         {\phi\left(\phibi\left(\frac{\alpha w_j}m\right)\right)}\\
\end{eqnarray*}
The $\underline w$ that solves these equations is given in (\ref{eq:wq1}).
Finally, solve for $c$ such that $\sum_i w_i = m$.
$\blacksquare$

\vspace{1cm}

{\sf Proof of Lemma \ref{lemma::discon}.}
Choose $K> 1$ such that
$1/(K+1) < 1/a - \epsilon$.
Choose $1> \gamma > (2\alpha-a)/K$.
Choose a small $c>0$.
Let
$\xi    = A + \sqrt{A^2 - 2c}$ and
$u     = B - \sqrt{B^2 - 2c}$
where
\begin{equation}
A     = \overline{\Phi}^{-1}\left(\frac{\alpha}{(m(\gamma K+a))}\right),\ \ \ \ 
B     = \overline{\Phi}^{-1}\left(\frac{K \alpha}{(m(\gamma K+a))}\right).
\end{equation}
Then $\rho(\xi)=1/a$ and
$\tilde{\rho}(\xi)=1/(K+1)$.
Now $d(Q,\tilde{Q})=\gamma$.
Taking $K$ sufficiently large and
$\gamma$ sufficiently close to $(2\alpha-a)/K$
makes $\gamma < \delta$. $\blacksquare$

\vspace{1cm}

{\sf Proof of Theorem \ref{thm::sparse-is-good}.}
The first statement follows easily by noting that the worst case
corresponds to choosing weight $B$ in the first term in $R(\xi)$ and
choosing weight $b$ in the second term in $R(\xi)$.
The rest follows by Taylor expanding
$R_{b,B}(\xi)$ around $b=1$. $\blacksquare$

\vspace{1cm}

{\sf Proof of Lemma \ref{lemma::safe}.}
With $b=0$, $R_{b,B}(\xi) \geq 0$ when
\begin{equation}\label{eq::bis0}
\overline{\Phi}(z_{B \alpha/m} - \xi) -2 \overline{\Phi}(z_{\alpha/m} - \xi) \geq 0.
\end{equation}
With $B\geq 2$, (\ref{eq::bis0}) holds at $\xi=0$.
The left hand side is increasing in $\xi$ for $\xi$ near 0 
but (\ref{eq::bis0}) does not hold at $\xi = z_{\alpha/m}$.
So (\ref{eq::bis0}) must holds in the interval
$[0,\xi_*]$.
Rewrite (\ref{eq::bis0}) as
$\overline{\Phi}(z_{B \alpha/m} - \xi) - \overline{\Phi}(z_{\alpha/m} - \xi) \geq \overline{\Phi}(z_{\alpha/m} - \xi)$.
We lower bound the left hand side and upper bound the right hand side.
The left hand side is
$\overline{\Phi}(z_{B \alpha/m} - \xi) - \overline{\Phi}(z_{\alpha/m} - \xi) =
\int_{ z_{B \alpha/m} - \xi}^{ z_{\alpha/m} - \xi} \phi(u) du \geq
(z_{\alpha/m} - z_{B\alpha/m})\phi(z_{\alpha/m}-\xi)$.
The right hand side
can be bounded using Mill's ratio:
$\overline{\Phi}(z_{\alpha/m} - \xi) \leq \phi(z_{\alpha/m} - \xi)/(z_{\alpha/m} - \xi)$.
Set the lower bound greater than the upper bound to obtain the stated result.
$\blacksquare$

\vspace{1cm}

It is convenient to prove
Theorem \ref{thm::drop-u} 
before proving Theorem
\ref{thm::restrict-u}.

{\sf Proof of Theorem \ref{thm::drop-u}}.
Let $c_*$ solve
\begin{equation}\label{eq::define-c-star}
\gamma\overline{\Phi}(\sqrt{2 c_*}) + a \overline{\Phi}\left( \frac{\xi}{2}+\frac{c_*}{\xi}\right)=
\frac{\alpha}{m}.
\end{equation}
We claim first that for any $c>c_*$,
there is no $u$ such that
the weights average to 1.
Fix $c > c_*$.
The weights average to 1 if and only if
\begin{equation}\label{eq::we-need}
\gamma\overline{\Phi}\left(\frac{c}{u}+\frac{u}{2}\right) +
a \overline{\Phi}\left( \frac{\xi}{2}+\frac{c}{\xi}\right)= \frac{\alpha}{m}.
\end{equation}
Since $c > c_*$ and since the second term
is decreasing in $c$, we must have
\begin{equation}
\overline{\Phi}\left(\frac{c}{u}+\frac{u}{2}\right)  > 
\overline{\Phi}(\sqrt{2 c_*}).
\end{equation}
The function 
$r(u)=\overline{\Phi}(c/u+u/2)$ is maximized at
$u=\sqrt{2c}$.
So $r(\sqrt{2c}) \geq r(u)$.
But $r(\sqrt{2c}) = \overline{\Phi}(\sqrt{2 c})$.
Hence
$\overline{\Phi}(\sqrt{2 c}) \geq r(u) \geq 
\overline{\Phi}(\sqrt{2 c_*}).$
This implies $c < c_*$ which is a contradiction.
This establishes that $\sup_u c(u) \leq c_*$.
On the other hand,
taking
$c=c_*$ and $u=\sqrt{2c_*}$ solves equation
(\ref{eq::we-need}). Thus $c_*$ is indeed the largest $c$ that solves the
equation which establishes the first claim.
The second claim follows by noting that
\begin{equation}
\gamma\overline{\Phi}(\sqrt{2 c_*}) + a \overline{\Phi}\left( \frac{\xi}{2}+\frac{c_*}{\xi}\right)=
\gamma\overline{\Phi}(\sqrt{2 c_*}) + O(a).
\end{equation}
Now set this expression equal to $\alpha/m$ and solve.
$\blacksquare$

\vspace{1cm}

{\sf Proof of Theorem \ref{thm::restrict-u}}.
Define $c_*$ as in (\ref{eq::define-c-star}).
If $u_*=\sqrt{2c_*} \leq \xi$ then the
the proof proceeds as in the previous proof.
So we first need to establish for which values of $\xi$
is this true.
Let $r(c) = \gamma\overline{\Phi}(\sqrt{2c}) + a\overline{\Phi}(\xi/2 + c/\xi)$.
We want to find out when 
the solution of $r(c)=\alpha/m$
is such that $\sqrt{2c} \leq \xi$, or
equivalently, $c \leq \xi^2/2$.
Now $r$ is decreasing in $c$.
Since $\gamma + a \geq \alpha/m$,
$r(-\infty) \geq \alpha/m$.
Hence there is a solution
with $c \leq \xi^2/2$
if and only if
$r(\xi^2/2) \leq \alpha/m$.
But
$r(\xi^2/2) = (\gamma + a)\overline{\Phi}(\xi)$
so we conclude that
there is such a solution if and only if
$(\gamma + a)\overline{\Phi}(\xi) \leq \alpha/m$, that is,
$\xi \geq z_{\alpha/(m(\gamma+a))}=\xi_0$.

Now suppose that
$\xi < \xi_0$.
We need to find $u\leq \xi$ to make $c$ as large as possible
in the equation
$v(u,c)\equiv\gamma \overline{\Phi}(u/2 + c/u) + a \overline{\Phi}(\xi/2 + c/\xi)=\alpha/m$.
Let $u_*=\xi$ and $c_*=\xi z_{\alpha/(m(\gamma+a))}-\xi^2/2$.
By direct substitution,
$v(u_*,c_*)=\alpha/m$ for this choice of $u$ and $c$
and clearly $u_*\leq \xi$ as required.
We claim that this is the largest possible $c_*$.
To see this, note that
$v(u,c)< v(u,c_*)$.
For $\xi\leq\xi_0$,
$v(u,c_*)$ is a decreasing function of $u$.
Hence,
$v(u,c) < v(u,c_*) \leq v(u_*,c_*)=\alpha/m$.
This contradicts the fact that
$v(u,c) =\alpha/m$.

For the second claim, note that the power
of the weighted test beats the power of Bonferroni
if and only if
the weight $w = (m/\alpha)\overline{\Phi}(\xi/2 + C(\xi)/2) \geq 1$
which is equivalent to
\begin{equation}\label{eq::this-must}
C(\xi) \leq \xi z_{\alpha/m} - \xi^2/2.
\end{equation}
When $\xi \leq \xi_0$,
$C(\xi) =\xi \xi_0 - \xi^2/2$.
By assumption, $\gamma + a \leq 1$ so that
$z_{\alpha/(m(\gamma+a))}\leq z_{\alpha/m}$ and
Now suppose that
$\xi_0 < \xi \leq \xi_*$.
Then $C(\xi)$ is the solution to
$r(c) = \gamma\overline{\Phi}(\sqrt{2c}) + a\overline{\Phi}(\xi/2 + c/\xi)=\alpha/m$.
We claim that
(\ref{eq::this-must}) still holds.
Suppose not.
Then, since $r(c)$ is decreasing in $c$,
$r(\xi z_{\alpha/m} - \xi^2/2) > r(C(\xi))=\alpha/m$.
But, by direct calculation, 
$r(\xi z_{\alpha/m} - \xi^2/2) > \alpha/m$
implies that $\xi > \xi_*$
which is a contradiction.
Thus (\ref{eq::thebonfclaim}) holds.

Finally, we turn to (\ref{eq::thenextclaim}).
In this case,
$C(\xi) = z^2_{\alpha/(m\gamma)}/2 + O(a)$.
The worst case power is
$\overline{\Phi}(C(\xi)/\xi - \xi/2)=\overline{\Phi}(z^2_{\alpha/(m\gamma)}/(2\xi) - \xi/2) + O(a)$.
The latter is increasing in $\xi$ and so is at least
$\overline{\Phi}(z^2_{\alpha/(m\gamma)}/(2\xi_*) - \xi_*/2) + O(a) =
\overline{\Phi}((z^2_{\alpha/(m\gamma)}/(2\xi_*) - \xi_*^2)/(2\xi_*))+O(a)$ as claimed.
The next two equations follow from standard tail approximations for
Gaussians. 
Specifically, a Gaussian quantile
$z_{\beta/m}$ can be written as
$z_{\beta/m}=\sqrt{2\log (m L_m/\beta)}$
where $L_m = c \log^a (m)$ for constants $a$ and $c$
(Donoho and Jin 2004).
Inserting this into the previous expression yields the final expression.
$\blacksquare$

\vspace{1cm}

{\sf Proof of Theorem \ref{thm::method1}.}
Setting $\pi(w,\xi_m) = \overline{\Phi}(\overline{\Phi}^{-1}(w\alpha/m) - \xi_m)$
equal to $1-\beta$ implies
$w = (m/\alpha)\overline{\Phi}(z_{1-\beta} + z_{\alpha/m})$
which is equal to $w_1$ as stated in the theorem.
The stated form of $w_0$ implies that the weights average to 1.
The stated solution thus satisfies
the restriction that a fraction $\epsilon$ have power at least $1-\beta$.
Increasing the weight of any hypothesis whose weight is $w_0$ necessitates
reducing the weight of another hypothesis.
This either reduces the minimum power of forces a hypothesis with power $1-\beta$
to fall below $1-\beta$.
Hence, the stated solution does in fact maximize the minimum power.
$\blacksquare$

\vspace{1cm}

The proof of Theorem \ref{thm::method2} is similar to the previous proof
and is omitted.

\vspace{1cm}

{\sf Proof of Theorem \ref{eq::newcontrol}.}
The familywise error is
\begin{eqnarray}
\mathbb{P}({\cal R}\cap {\cal H}_0) & \leq &
\sum_{j\in {\cal H}_0}\mathbb{P}\left(P_j\left(T_j^{(1)},T_j^{(2)}\right)\leq \frac{w_j(T_j^{(1)})\alpha}{m}\right)\\
&=& \sum_{j\in {\cal H}_0}
\mathbb{E}
\left(\mathbb{P}\left(
P_j\left(t_j^{(1)},T_j^{(2)}\right)\leq \frac{w_j(t_j^{(1)})\alpha}{m}\Biggm| T_j^{(1)}=t_j^{(1)}\right)\right)\\
&=& \sum_{j\in {\cal H}_0}
 \mathbb{E}\left(\mathbb{P}
  \left(\overline{\Phi}(t_j^{(1)},T_j^{(2)})\leq 
 \frac{w_j(t_j^{(1)})\alpha}{m}\Biggm| T_j^{(1)} =t_j^{(1)}\right)\right)\\
&=& \sum_{j\in {\cal H}_0}
 \mathbb{E}\left(\mathbb{P}
  \left((1-b)^{1/2}T_j^{(2)} \geq \overline{\Phi}^{-1}\left(\frac{w_j(t_j^{(1)})\alpha}{m}\right)-b^{1/2} t_j^{(1)}\right)\right)\\
&=& \sum_{j\in {\cal H}_0}
 \mathbb{E}\left(\overline{\Phi}
 \left(\frac{\overline{\Phi}^{-1}\left(\frac{w_j(t_j^{(1)})\alpha}{m}\right)-b^{1/2} t_j^{(1)}}{(1-b)^{1/2}}\right)
\right)\\
&=& \sum_{j\in {\cal H}_0}
 \mathbb{E}\left(\overline{\Phi}
 \left(\frac{ \frac{\hat\xi_j}{2} + \frac{c(t^{(1)})}{\hat\xi_j} - b^{1/2} t_j^{(1)}}{(1-b)^{1/2}}\right)\right)\\
& \leq & \sum_{j=1}^m
 \mathbb{E}\left(\overline{\Phi}
 \left(\frac{ \frac{\hat\xi_j}{2} + \frac{c(t^{(1)})}{\hat\xi_j} - b^{1/2} t_j^{(1)}}{(1-b)^{1/2}}\right)\right)\\
& \leq & \mathbb{E} \sum_{j=1}^m 
 \left(\overline{\Phi}
 \left(\frac{ \frac{\hat\xi_j}{2} + \frac{c(t^{(1)})}{\hat\xi_j} - b^{1/2} t_j^{(1)}}{(1-b)^{1/2}}\right)\right)\\
&=& \alpha.\ \ \ \blacksquare
\end{eqnarray}

\newpage

\section*{References}

\begingroup
\parindent=0em
\baselineskip=17pt
\parskip=6pt
\everypar={\hangindent=1em\hangafter=1\noindent}

Benjamini, Y. and Hochberg, Y. (1995).
Controlling the false discovery rate: 
A practical and powerful approach to multiple testing.
{\em Journal of the Royal Statistical Society, Series B}, 57, 289--300.

Benjamini, Y.~and Hochberg, Y. (1997). 
Multiple Hypothesis Testing with Weights, 
{\em Scandinavian Journal of Statistics}, {\bf 24}, 407--418.

Chen, James J. and Lin, Karl K. and Huque, Mohammad and Arani, Ramin B. (2000).
Weighted $p$-value adjustments for animal carcinogenicity trend test.
{\em Biometrics}, 56, 586--592.

Donoho, D. and Jin, J. (2004).
Higher criticism for detecting sparse heterogeneous mixtures.
{\em The Annals of Statistics}, {\bf 32}, 962--994.

Genovese, C.~R., Roeder, K. ~and Wasserman, L.~(2005).
False Discovery Control with P-Value Weighting.
To appear: {\em Biometrika}.

Holm S.~(1979). A simple sequentially rejective multiple test procedure. 
{\em Scandinavian Journal of Statistics}, {\bf 6}, 65--70.

Roeder, Bacanu, Wasserman and Devlin (2005).
Using Linkage Genome Scans to
Improve Power of Association in Genome Scans.
{\em The American Journal of Human Genetics.}
{\bf 78}.

Rubin, D, van der Laan, M. and Dudoit, S. (2005). 
Multiple testing procedures which are optimal at a simple alternative. 
Technical report 171,
Division of Biostatistics,
School of Public Health,
University of California, Berkeley .

Skol A.D., Scott L.J., Abecasis G.R. and Boehnke M. (2006) Joint
analysis is more efficient than replication-based analysis for
two-stage genome-wide association studies.  {\em Nat Genet.}
38:390-394.

Signoravitch, J. (2006).
Optimal multiple testing under the general linear model.
Technical report.
Harvard Biostatistics.

Storey J.D. (2005).  The optimal discovery procedure: A new approach
to simultaneous significance testing.  UW Biostatistics Working Paper
Series, Working Paper 259.

Thomas D.C., Haile R.W. and Duggan D. (2005)       
Recent developments in genomewide association scans: 
a workshop summary and review.
{\em Am J Hum Genet}  77:337-345

\endgroup

\end{document}